\documentclass{proc-l}
\usepackage{amsfonts}
\usepackage{amssymb}
\normalsize 
 \divide\textheight by
\baselineskip \multiply\textheight by \baselineskip
\newcommand{\HHH}{{\mathcal H}ol}

\newcommand{\C}{\mathbb{C}}

\newcommand{\te}{Teich\-m\"{u}l\-ler}
\newtheorem{thm}{Theorem}[section]
\newtheorem{lemma}[thm]{Lemma}

\theoremstyle{definition}
\newtheorem{definition}[thm]{Definition}
\newtheorem{example}[thm]{Example}

\theoremstyle{remark}

\numberwithin{equation}{section}

\newtheorem{cor}{Corollary}[section]

\begin{document}
\title[Random Iterations]{Random holomorphic iterations and degenerate subdomains of the unit disk}
\author{Linda Keen}
\address{Department of Mathematics\\Lehman College and Graduate Center, CUNY\\ Bronx NY 10468}
\email{linda@lehman.cuny.edu}
\thanks{Paritally supported by a PSC-CUNY Grant}
\author{Nikola Lakic}
\address{Department of Mathematics\\Lehman College and Graduate Center, CUNY\\ Bronx NY 10468}
\email{nlakic@lehman.cuny.edu}
\thanks{Partially supported by NSF grant DMS
0200733 } 
\date{\today}

\begin{abstract} Given a random sequence of
holomorphic  maps $f_1,f_2,f_3,\ldots$ of the unit disk $\Delta$
to a subdomain $X$, we consider the compositions
$$F_n=f_1 \circ f_{2} \circ \ldots f_{n-1} \circ f_n.$$
 The sequence $\{F_n\}$ is  called the {\em iterated function
system} coming from the sequence $f_1,f_2,f_3,\ldots.$ We prove
that a sufficient condition on the domain $X$ for all limit
functions of any $\{F_n\}$ to be constant is also necessary. We
prove the condition is a quasiconformal invariant.  Finally, we
address the question of uniqueness of limit functions.
\end{abstract}

\maketitle

\begin{section}{Introduction} \label{sec:introduction}
 Suppose that we are given a random sequence of
holomorphic self maps $f_1,f_2,f_3,\ldots$ of the unit disk
$\Delta$.  We consider the compositions
$$F_n=f_1 \circ f_{2} \circ \ldots f_{n-1} \circ f_n.$$
 The sequence $\{F_n\}$ is  called the {\em iterated function
system} coming from the sequence $f_1,f_2,f_3,\ldots.$    By
Montel's theorem (see for example \cite{CG}), the sequence $F_n$
is a normal family, and every convergent subsequence converges
uniformly on compact subsets of $\Delta$ to a holomorphic function
$F$.  The limit functions $F$ are called accumulation points.
Therefore every accumulation point is either an open self map of
$\Delta$ or a constant map. The constant accumulation points may
be located either inside $\Delta$ or on its boundary.

We may look at the iterated function system as a dynamical system
acting on $\Delta$. If $z$ is an arbitrary point of $\Delta$, its
orbit under the iterated function system, $F_n(z)$, has $F(z)$ as
an accumulation point. Hence, if the only limit functions are
constants, the orbits of all points tend to periodic cycles. As we
will see, we can find conditions so that whether this happens
depends only on the subdomain $X \subset \Delta$ and not on the
particular system chosen from $\HHH(\Delta,X)$.

If all maps in the iterated system are the same, the well known
Denjoy-Wolff Theorem determines all possible accumulation points.

\vspace{.2in}

\noindent {\bf The Denjoy-Wolff Theorem.} \ \ {\em Let $f$ be a
holomorphic self map of the unit disk $\Delta$ that is not a
conformal automorphism.  Then the iterates $f^{\circ n}$ of $f$
converge locally uniformly in $\Delta$ to a constant value $t,$
where $|t| \leq 1.$}

\vspace{.2in}

Therefore, whenever $f$ is not a biholomorphic isometry of
$\Delta$, the system $f^{\circ n}$ has a constant accumulation
point.

Many articles have studied possible generalizations of the
Denjoy-Wolff Theorem to  iterated function systems. One result of
Lorentzen  and Gill  is

\vspace{.2in}

\noindent {\bf Theorem(GL})( \cite{Gi},\cite{Lo}) \ {\em If an
iterated function system is formed from functions in
$\HHH(\Delta,X)$ where $X$ is relatively compact in $\Delta$ then
the system $F_n$ converges locally uniformly in $\Delta$ to a
unique constant. This constant is, of course, located in the
relatively compact set $X$.}

\vspace{.2in}

We say that a subdomain $X$ of $\Delta$ {\em is degenerate} if
every iterated function system generated by a sequence of maps
from $\Delta$ to $X$ has only constant accumulation points.
Therefore, any relatively compact subdomain of $\Delta$ is
degenerate, and moreover, each system has a  constant limit.

A recent study of iterated function systems (see \cite{BCMN})
introduced new degenerate subdomains that are not relatively
compact in $\Delta$.  It also considered more general iterated
function systems formed from maps in $\HHH(\Omega,X)$ where $X
\subset \Omega $ are arbitrary plane domains such that $\Omega$
(and hence also $X$) admits $\Delta$ as a universal cover. We call
such  domains {\em hyperbolic domains}. In this context we say
{\em $X$ is degenerate in $\Omega$ if all the accumulation points
of every iterated function system coming from $\HHH(\Omega,X)$ are
constant}. If $X$ is relatively compact in $\Omega$, we can  apply
theorem GL to the universal covers of $\Omega$ and $X$ to show
that $X$ is degenerate; again the interesting case is when $X$ is
not relatively compact in $\Omega$.

The study cited above uses the Poincar\'e metrics on $\Delta$,
$\Omega$ and $X$, denoted respectively by $\rho=\rho_{\Delta},
\rho_{\Omega}$ and $\rho_{X}$, as an important tool. The authors
extended the classical Euclidean notion of a Bloch subdomain of
the Euclidean plane to the hyperbolic setting. A
 subdomain $X$ of $\C$ is a Euclidean Bloch subdomain if, and only if,
there is an upper bound on the radii of the disks lying in $X.$ In
the hyperbolic context we define
\begin{definition} Let $R(X,\Omega)$ be the supremum of all radii (measured with
respect to $\rho_{\Omega}$) of hyperbolic subdisks of $\Omega$
that are contained in $X.$ A subdomain $X$ of $\Omega$ is called
{\em a $\rho$-Bloch subdomain of $\Omega$} if
$$R(X,\Omega)<\infty.$$ $R(X,\Omega)$ is called the {\em
$\rho$-Bloch radius of $X$ in $\Omega$}.
\end{definition}

In their paper, Beardon, Carne, Minda and Ng  prove

\vspace{.2in}

\noindent {\bf Theorem BCMN.} \ \  {\em If $X$ is a $\rho$-Bloch
subdomain of $\Omega,$ then $X$ is degenerate in $\Omega.$}

\begin{example}
{\rm Suppose that $\Omega$ is a plane domain obtained by removing
at least two but only finitely many points from the whole complex
plane. Suppose that $X$ is obtained by removing at least one but
finitely many points from $\Omega.$ Then any holomorphic map from
$\Omega$ to $X$ has only removable singularities and therefore
extends to a rational map of the whole complex sphere. It is easy
to see that such a map can not have an image in a strictly smaller
subdomain, unless it is a constant map. Therefore, every
holomorphic map from $\Omega$ to $X$ is a constant map and $X$ is
degenerate in $\Omega.$ Small punctured disks about each of the
complementary points of $\Omega$ in $\C$ contain (non-schlicht)
hyperbolic disks of arbitrarily large radii so $X$ is a
non-$\rho$-Bloch subdomain of $\Omega.$ This example shows that
the generalization of the converse of theorem~BCMN does not hold:
$X$ is non-$\rho$-Bloch in $\Omega$ but is degenerate.
Furthermore, any iterated function system $F_n$ converges to the
constant $f_1(z)$ and therefore the converse of the generalized
version of Theorem~GL does not hold either. }
\end{example}

The situation is better, however, in the case when $\Omega$ is
simply connected.   In section~\ref{sec:nonconstantlimits} we
prove
\begin{thm} \label{thm:theorem7}
Suppose that $X \subset \Delta$ is not a $\rho$-Bloch subdomain of
$\Delta.$  Then $X$ is not degenerate in $\Delta$.
\end{thm}
Therefore, theorem \ref{thm:theorem7} together with theorem~BCMN
implies

\begin{cor} \label{coro}
$X$ is a $\rho$-Bloch subdomain of $\Delta$ if, and only if, $X$
is degenerate in $\Delta$.
\end{cor}

In section~\ref{qcinvariant} we show that $\rho$-Bloch subdomains
are quasiinvariant with respect to the unit disk.  More precisely
we prove

\begin{thm} \label{qcthm}
If $f$ is a quasiconformal self homeomorphism of the unit disk
$\Delta,$ then $f$ maps every $\rho$-Bloch subdomain of $\Delta$
onto a $\rho$-Bloch subdomain of $\Delta.$
\end{thm}

This theorem together with corollary \ref{coro} implies

\begin{cor} \label{coro2} If $f$ is a quasiconformal self
homeomorphism of the unit disk $\Delta,$ then $f$ maps degenerate
subdomains of $\Delta$ onto degenerate subdomains of $\Delta.$
\end{cor}

In section~\ref{sec:nonuniquelimits} we turn our attention to the
question of the uniqueness of limits of iterated function systems.
We prove the converse of theorem GL:
\begin{thm} \label{thm:theorem8}
Suppose that $X$ is any subdomain of the unit disk $\Delta$ that
is not relatively compact in $\Delta.$ Then there exists a
sequence $f_n$ of holomorphic mappings from $\Delta$ to $X$ such
that the iterated system $F_n=f_1\circ\ldots\circ f_n$ has more
than one accumulation point.
\end{thm}

 Several recent articles are concerned with the study of
iterated function systems
 and its applications, see for example, \cite{BCMN},\cite{PU},\cite{MU},\cite{SU}.
We would like to thank Fred Gardiner for numerous helpful
discussions and remarks on an earlier version of this paper and
Jonathan Brezin for his editorial advice.

\end{section}

\begin{section}{Non-constant accumulation points}
\label{sec:nonconstantlimits} In this section we prove
theorem~\ref{thm:theorem7}. We begin by proving two preparatory
lemmas.  The first gives an estimate depending on the $\rho$-Bloch
radius that relates the distances between relatively close points
in the $\rho=\rho_{\Delta}$ and $\rho_X$ metrics.  We  normalize the density function
for $\rho$ by $\rho(z)=\frac{1}{1-|z|^2}$.
\begin{lemma}
\label{lemma:prep1}  Let $a$ be a point in a subdomain $X$ of
$\Delta.$ Let $C=R(X,\Delta,a)$ be the radius of the largest
$\rho$-disk with center at $a$ which is inside $X.$   If $z$ is
another point in $X,$ with
$$\rho(a,z)<1<C,$$ then
$$\rho_{X}(a,z)\leq (1+\epsilon)\rho(a,z)$$
where $$\epsilon=\epsilon(C)\rightarrow 0 \mbox{  as  }C
\rightarrow \infty.$$
\end{lemma}
\begin{proof}
By applying the M\"obius transformation
$$A(z)=\frac{z-a}{1-\overline{a}z},$$
 we may assume that $a=0.$  Suppose that $C=R(X,\Delta,0)>1>\rho(0,z)$
 for some point $z$ in $X.$  Let $D$ be a disk in $\Delta$ with
 center at $0$ and  $\rho$-radius $C.$  Then $D \subset X,$  so that
   $\rho_{X} \leq \rho_D.$  Therefore, an easy calculation shows
$$\rho_{X}(0,z) \leq \rho_D(0,z)= \int_0^{|z|} \frac{c}{c^2-t^2} dt$$
where $c$ is the Euclidean radius of $D.$  Obviously $c
\rightarrow 1$ as $C \rightarrow \infty.$ Therefore,
$$\rho_{X}(0,z)\leq \int_0^{|z|} \frac{1}{1-t^2}
(c+\frac{c(1-c^2)}{c^2-t^2})dt\leq$$
$$ \int_0^{|z|} \frac{1}{1-t^2} (c+\frac{c(1-c^2)}{c^2-|z|^2})dt
\rightarrow \int_0^{|z|} \frac{1}{1-t^2}  dt= \rho(0,z)$$ as $c
\rightarrow 1,$ and the lemma follows.
\end{proof}
The second preparatory lemma is about the contraction properties
of  Blaschke product maps of degree two.
\begin{lemma}
\label{lemma:prep2} Let $c \neq 0$ be any point in  $\Delta$ such
that $\rho(0,c) < 1$.  If
$$A_a(z)=\frac{z(z-a)}{1-\overline{a}z},$$ then $A^{-1}(c)=\{z_1,z_2\}$
 (that is, $A_a(z_1)=A_a(z_2)=c$)
and
$$\rho(0,z_1)=\rho(a,z_2)\rightarrow \rho(0,c) \mbox{  as  } |a| \rightarrow 1.$$
\end{lemma}
\begin{proof}
The two points $z_1$ and $z_2$ are the solutions of the equation
$A_a(z)=c$ and thus roots of
\begin{equation}z^2-z(a-\overline{a}c)-c=0
\label{eqtn}
\end{equation}
Therefore, we have
\begin{equation}
\label{c5}  z_1z_2=-c
\end{equation}
Since $A_a(z_2)=c,$ by (\ref{c5}) we have
$$\frac{z_2-a}{1-\overline{a}z_2}=\frac{c}{z_2}=-z_1,$$
Thus,
$$\rho(a,z_2)=\rho(0,-z_1)=\rho(0,z_1).$$
Solving  equation (\ref{eqtn}), we obtain
$$2z_{1,2}=a-\overline{a}c\pm \sqrt{a^2+\overline{a}^2c^2+2c(2-|a|^2)}.$$
Hence,
$$2az_{1,2}=a^2-|a|^2c\pm\sqrt{a^4+|a|^4c^2+2a^2c(2-|a|^2)},$$
$$\lim_{|a|\rightarrow 1}|2z_{1,2}|=\lim_{|a|\rightarrow 1} |a^2-c\pm\sqrt{a^4+c^2+2a^2c}|.$$
Therefore we may order $z_1$ and $z_2$ so that
$$\lim_{|a|\rightarrow 1}|2z_1|= |2c|\quad \mbox{ and }\quad
 \lim_{|a|\rightarrow 1}|2z_2|= \lim_{|a|\rightarrow 1}|2a^2|= 2$$
Finally, $\lim_{|a|\rightarrow 1} \rho(0,z_1)=\rho(0,c),$
finishing the proof of this lemma.
\end{proof}

\vspace{.2in}

\noindent {\bf  Theorem~\ref{thm:theorem7} } \ {\em Suppose that
$X \subset \Delta$ is not a $\rho$-Bloch subdomain of $\Delta.$
Then $X$ is not degenerate in $\Delta$.}

\begin{proof}
 Let $X$ be any non-$\rho$-Bloch domain.  We are going to construct
an iterated function system from $\HHH(\Delta,X)$ with  a
nonconstant accumulation point. Pick any two distinct points $a_0$
and $w_0$ in $X$ such that $\rho_{X}(a_0,w_0)<1/2.$ We will
recursively find functions $f_n$ such that  the iterated function
system $F_n$ will have  limit function $F$ that satisfies
$F(0)=a_0$ and $F(\tilde{w})=w_0$ for some $\tilde{w}$,
$\rho(0,\tilde{w})<1$.

First let $\pi_1$ be a universal covering map from  $\Delta$ onto
$X$ such that $\pi_1(0)=a_0.$ Then  there exists a point $c_0 \in
\Delta,$ such that $\pi_1(c_0)=w_0$ and
\begin{equation}
\rho(0,c_0)=\rho_{X}(a_0,w_0) \label{aaa}
\end{equation}
For any  choice of $a_1 \in X$, the Blaschke product $A_{a_1}$
produces two points $w_1$ and $\tilde{w}_1$ in $\Delta$ as the
preimages of $c_0:$
\begin{equation}
\label{a1} A_{a_1}(\tilde{w_1})=A_{a_1}(w_1)=c_0 \mbox{ and
}\rho(0,\tilde{w_1})=\rho(a_1,w_1)
\end{equation}
Define $f_1=\pi_1 \circ A_{a_1}$ so that
\begin{equation}
\label{b1b}
  f_1(0)=f_1(a_1)=a_0,
\end{equation}

\begin{equation}
\label{b1bb} f_1(w_1)=f_1(\tilde{w}_1)=w_0. \end{equation}

%
%
%

We need to make sure that the point $w_1$ belongs to $X$. To do
this we use the preparatory lemmas.

Let $\epsilon_n \rightarrow 0$ be a sequence such that
$\Pi_1^{\infty}(1+\epsilon_n)^2 \leq 2$.  Since $X$ is
non-$\rho$-Bloch, by lemma~\ref{lemma:prep2},  we can choose $a_1
\in X$ so that
 $|a_1|$ is close enough to $1$
so that
\begin{equation}
\rho(0,\tilde{w_1})=\rho(a_1,w_1)<(1+\epsilon_1)
\rho(0,c_0)=(1+\epsilon_1)\rho_X(a_0,w_0) \label{a1a}
\end{equation}
Moreover, we may  assume that $R(X,\Delta,a_1)$ is greater
than $1$ so that  formula~(\ref{a1a}) implies that $w_1 \in X$ and also 
that it is large enough so that by lemma~\ref{lemma:prep1} and
(\ref{a1a}) we get
\begin{equation}\label{a1aaa}
\rho_X(a_1,w_1) < (1+\epsilon_1)\rho(a_1,w_1) < (1+\epsilon_1)^2
\rho_X(a_0,w_0)<1
\end{equation}
Then, inductively, by our choice of $\epsilon_n$, there exist
points $a_n,w_n \in X$ and $\tilde{w}_n \in \Delta$ such that

\begin{equation}
\label{d1d}  f_n(0)=f_n(a_n)=a_{n-1} \quad \mbox{ and }\quad
 f_n(w_n)=f_n(\tilde{w}_n)=w_{n-1},
\end{equation}
 
   \begin{equation} \label{d1ddd} \rho(a_n,w_n)=\rho(0,\tilde{w}_n) <
(1+\epsilon_n)\rho_X(a_{n-1},w_{n-1})
\end{equation}
and
\begin{equation}
\label{d1dddd}\rho_X(a_n,w_n) <(1+\epsilon_n)
  \rho(a_{n},w_{n})<(1+\epsilon_n)^2\rho_X(a_{n-1},w_{n-1})
\end{equation}
Therefore
\begin{equation}\label{eqn:summ1}
\rho(0,\tilde{w}_n) < \Pi_1^n(1+\epsilon_i)\rho(0,c_0)<1
\end{equation}
\begin{equation} \label{eqn:summ2}
\rho_{X}(a_n,w_n)< \Pi_1^n(1+\epsilon_i)^2\rho(0,c_0)<1
\end{equation}
Now if $F_n=f_1\circ f_2 \circ f_3 \circ \ldots f_n,$ the
formulas~{\ref{d1d}} yield

\begin{equation}
\label{wn}
F_n(0)=a_0 \quad \mbox{ and } \quad F_n(\tilde{w}_n)=w_0, 
\end{equation}

By Montel's theorem, ${F_n}$ is a normal family so that a
subsequence of $F_n$ converges uniformly on compact subsets of
$\Delta$ to a holomorphic limit function $F.$ Therefore,
equations~(\ref{wn}) yield 
$$F(0)=a_0 \quad \mbox{ and }\quad F(\tilde{w})=w_0,$$
where $\tilde{w}$ is an accumulation point of the sequence
$\tilde{w}_n.$ Since, by equation~\ref{eqn:summ1},
$\rho(0,\tilde{w}_n)<1$ for all $n,$ the point $\tilde{w}$ belongs
to $\Delta.$ This implies that $F$ is a nonconstant function.
\end{proof}
\end{section}

\begin{section}{Degeneracy in the unit disk is a quasiconformal invariant}  \label{qcinvariant}

Since conformal homeomorphisms  of the unit disk  
 are hyperbolic isometries, conjugating an iterated function system  by one preserves degeneracy.  Here we show that the same is true for quasiconformal homeomorphisms. 
\vspace{.2in}

\noindent {\bf  Theorem~\ref{qcthm} } \ {\em If $f$ is a
quasiconformal self homeomorphism of the unit disk $\Delta,$ then
$f$ maps every $\rho$-Bloch subdomain of $\Delta$ onto a
$\rho$-Bloch subdomain of $\Delta.$}

\begin{proof}
Suppose that $f$ is a $K-$qusiconformal map from the unit disk
onto the unit disk and let $X$ be a non-$\rho$-Bloch subdomain of
the unit disk.  Then for every positive integer $n$ there exists a
point $p_n \in X$ and a hyperbolic disk $D_n \subset X$ centered
at $p_n$ with $\rho$-radius $n$.  Suppose that $f(X)$ is
$\rho$-Bloch.\; then $R(f(X),\Delta) < \infty.$  Let $F_n= h_n
\circ f \circ g_n,$ where
$$g_n(z)= \frac{z+p_n}{1+\overline{p_n}z} \quad \mbox{ and }\quad
 h_n(z)= \frac{z-f(p_n)}{1-\overline{f(p_n)}z}.$$ Then $F_n$ is a
K-quasiconformal map from the unit disk onto itself and $F_n(0)=0.$

 Since $h_n$ is an isometry for $\rho$, there exists a point $q_n$ in $\Delta \setminus
h_n(f(X))$ such that $\rho(q_n,0)\leq R(f(X),\Delta).$ Since $g_n$
is an isometry for $\rho$, the hyperbolic disk $g_n^{-1}(D_n)$
with center at 0 and radius $n$ is contained in $g_n^{-1}(X).$
 Therefore, 
$F^{-1}_n(q_n)$ must be outside the hyperbolic disk $g^{-1}(D_n).$
Thus, $|F^{-1}_n(q_n)| \rightarrow 1$  while $q_n$ stays bounded
inside the unit disk. This contradicts the fact that the family of
K-quasiconformal self homeomorphisms of the unit disk fixing 0 is a
normal family (see for example (\cite{QC}).
\end{proof}
Theorem (\ref{qcthm}) together with corollary (\ref{coro})
obviously implies

\vspace{.2in}

\noindent {\bf  Corollary~\ref{coro2} } \ {\em If $f$ is a
quasiconformal self homeomorphism of the unit disk $\Delta,$ then
$f$ maps degenerate subdomains of $\Delta$ onto degenerate
subdomains of $\Delta.$}

\end{section}

\begin{section}{Non-Uniqueness of limit points}
\label{sec:nonuniquelimits}

In this section we study the question of the uniqueness of the
limit points of an iterated function system  in $\HHH(\Delta,X)$.
A subdomain $X$ of $\Delta$ is degenerate if and only if it is a
$\rho$-Bloch domain. We show that there are sequences in
$\HHH(\Delta,X)$ that have more than one accumulation point. By
theorem~GL, if the subdomain $X$ is
 relatively compact the limit points are unique, so we assume that $X$ is an
arbitrary non-relatively compact  subdomain of $\Delta.$

\vspace{.1in}

\noindent {\bf Theorem~\ref{thm:theorem8}} \ {\em Suppose that $X$
is any subdomain of the unit disk $\Delta$ that is not relatively
compact in $\Delta.$ Then there exists a sequence $f_n$ of
holomorphic mappings from $\Delta$ to $X$ such that  the iterated
system $F_n=f_1\circ\ldots\circ f_n$ has more than one
accumulation point.}

\begin{proof}
Let $X$ be any  subdomain of  $\Delta$ that is not relatively
compact.  In our construction, all maps $f_i$ will be different
universal covering maps from $\Delta$ onto $X.$ We start with an
arbitrary point $a$ in $X.$  If we choose a point $a_1 \neq a \in
X$, we can find $f_1$ such that $f_1(a)=a_1$.  Then $f_1$ is
defined up to a (hyperbolic) rotation about $a$.

Let $f_2$ be a covering map from $\Delta$ onto $X$ such that
$f_2(a)=a_1$ and such that there is an $a_2 \in X$ with
$f_2(a_2)=a$. There is such an $a_2$ because $X$ is not relatively
compact in $\Delta$ so every hyperbolic circle with center at $a$
intersects $X$ and any covering map sending $a$ to $a_1$ is
defined only up to rotation about $a$. Since covering maps are
local isometries we may also assume $\rho(a_2,a)=\rho_{X}(a,a_1).$

Continuing this process we obtain a sequence of covering maps
$f_n$ and a sequence of points $a_n$ in $X$ such that

\begin{equation}
f_n(a)=a_{n-1} \mbox{ and } f_n(a_n)=a  \label{aan}
\end{equation}
for all $n.$  We study the even and  odd subsequences of the
iterated function system $F_n=f_1\circ\ldots\circ f_n$.   The
equalities (\ref{aan}) imply that the even subsequence $F_{2n}$
satisfies $F_{2n}(a)=a,$ and the odd subsequence $F_{2n+1}$
satisfies $F_{2n+1}(a)=a_1\neq a.$ Therefore these two
subsequences have different accumulation points.
 \end{proof}

Theorems~GL,BCMN,~\ref{thm:theorem7},~\ref{qcthm}
and~\ref{thm:theorem8} immediately yield the following corollary.
\begin{cor}
If $X$ be a subdomain of the unit disk $\Delta$, then:\\

\noindent (1) \ All accumulation points of any iterated function
system of maps in $\HHH(\Delta,X)$ are constant functions if, and
only if, $X$ is a $\rho$-Bloch-subdomain of $\Delta$.\\

\noindent (2) \ These accumulation points are unique if, and only
if, $X$ is a relatively compact subdomain of $\Delta.$\\

\noindent (3) \ The properties in both (1) and (2) are preserved
under quasiconformal self homeomorphisms of the unit disk.
\end{cor}
\end{section}

\vspace{.2in}

\bibliography{newver}

\begin{thebibliography}{10}
\bibitem{Ahl}L.~V.~Ahlfors
\newblock \underline{Complex Analysis},
\newblock  McGrawHill, (1953)

\bibitem{QC}L.~V.~Ahlfors
\newblock \underline{Lectures on Quasiconformal Mappings},
\newblock  Van Nostrand, (1966)

\bibitem{BCMN} A.~F.~Beardon, T.~K.~Carne, D.~Minda and T.~W.~Ng,
\newblock Random iteration of analytic maps,
\newblock preprint.

\bibitem{CG} L.~Carleson and T.~W.~Gamelin,
\newblock \underline{Complex Dynamics},
\newblock Springer-Verlag (1993).


\bibitem{G} F.~P.~Gardiner, \newblock oral communication.


\bibitem{G2} \underline{\hspace{.30in}},
\newblock \underline{\te\ Theory and Quadratic Differentials,}
\newblock Wiley-Interscience, 1987.


\bibitem{GL} F.~P.~Gardiner and N.~Lakic,
\newblock  \underline{Quasiconformal \te\ Theory,}
\newblock AMS Mathematical Surveys and Monographs,
{\bf 76}, 2000.

\bibitem{GLii} F.~P.~Gardiner and N.~Lakic,
\newblock  Comparing Poincar\'e distances
\newblock Annals of Math., {\bf 154}, 2001, 245--267
\bibitem{Gi} J.~Gill,
\newblock Compositions of analytic functions of the form $F_n(z)=F_{n-1}(f_n(z)), \ f_n(z)
\rightarrow f(z),$
\newblock  J. Comput. Appl. Math., {\bf 23} (2), 1988, 179--184

\bibitem{KLI} L.~Keen and N.~Lakic
\newblock Forward Iterated Function Systems
\newblock To appear, Proc. 2002 Workshop on Complex Dynamics,
Morningside Institute, Beijing China

\bibitem{KLbook} L.~Keen and N.~Lakic
\newblock \underline{An Introduction to Hyperbolic Geometry in the
Small}
\newblock In preparation

\bibitem{Lo} L.~Lorentzen,
\newblock Compositions of contractions,
\newblock  J. Comput. Appl. Math., {\bf 32}  1990, 169--178




\bibitem{PU} D.~Mauldin, F.~Przytycki and M.~Urbanski,
\newblock Rigidity of conformal iterated function systems,
\newblock  Compositio Math, {\bf 129}  2001, 273--299

\bibitem{MU} V.~Mayer, D.~Mauldin and M.~Urbanski,
\newblock Rigidity of connected limit sets of iterated function systems,
\newblock  Mich. Math J., {\bf 49}  2001, 451--458

\bibitem{SU} B.~Solomyak and M.~Urbanski,
\newblock $L^q$ densities for measures associated with parabolic iterated function systems with
overlaps,
\newblock  Indiana J. Math., {\bf 50}  2001

\bibitem{SV} T.~Sugawa and M. Vourinen,
\newblock Some inequalities for the Poincar\'e metric of plane domains,
\newblock  preprint


\end{thebibliography}

\def\noopsort#1{} \def\printfirst#1#2{#1} \def\singleletter#1{#1}
\def\switchargs#1#2{#2#1} \def\bibsameauth{\leavevmode\vrule height .1ex
depth 0pt width 2.3em\relax\,}

\end{document}